\documentclass[11pt]{amsart}

\usepackage{amssymb} % mathematical fonts (in amsfonts)
\usepackage[mathscr]{eucal} % load euler fonts (in amsfonts)
\usepackage[matrix,arrow,tips,curve,ps]{xy}
\usepackage{amsmath} %this is needed for \rxar (\xrightarrow)
\usepackage{epsfig}
\usepackage{amscd}
\usepackage[usenames,dvipsnames,svgnames,table]{xcolor}
\usepackage{tikz}
\usepackage{tikz-cd}
\usepackage{mathrsfs}
\usepackage{verbatim}

\newtheorem{TEO}{Theorem}[section]
\newtheorem{PROP}[TEO]{Proposition}
\newtheorem{LEM}[TEO]{Lemma}

\theoremstyle{definition}
\newtheorem{remark}[TEO]{Remark}
\newtheoremstyle{dico}% name of the style to be used
 {\baselineskip}   % ABOVESPACE
  {\topsep}   % BELOWSPACE
  {}  % BODYFONT
  {0pt}       % INDENT (empty value is the same as 0pt)
  {} % HEADFONT
  {.}         % HEADPUNCT
  {5pt plus 1pt minus 1pt} % HEADSPACE
  {}          % CUSTOM-HEAD-SPEC
\theoremstyle{dico}

\numberwithin{equation}{section}

\def\OO{{\mathcal O}}

\newcommand\dual{\mathrel{\raise3pt\hbox{$\underline{\mathrm{\thinspace d
\thinspace}}$}}}

\newcommand\proj{\mathbb P}
\newcommand\Z{\mathbb Z}

\newcommand\Co{\mathbb C}

\newcommand\rank{\operatorname{rank}}

\newcommand{\Nm}{\operatorname{Nm}}

\newcommand{\PP}{{\mathbb P}}

\renewcommand{\phi}{\varphi}

\newcommand{\ra}{\rightarrow}

\newcommand{\tC} {{\tilde{C} }}

\newcommand{\mihi}[1]{}

\begin{document}

\footnote{
The first  author was partially supported by MIUR PRIN 2015``Geometry of Algebraic Varieties''.
   The second author was partially supported by MIUR PRIN
   2015 ``Moduli spaces and Lie theory''.   The authors were also
   partially supported by GNSAGA of INdAM.
 }

\title{Second fundamental form of the Prym map in the ramified case}

\author[E. Colombo]{Elisabetta Colombo}
\address{Dipartimento di Matematica,
Universit\`a di Milano, via Saldini 50,
     I-20133, Milano, Italy } \email{{\tt
elisabetta.colombo@unimi.it}}

\author[P. Frediani]{Paola Frediani}
\address{ Dipartimento di Matematica, Universit\`a di Pavia,
via Ferrata 5, I-27100 Pavia, Italy } \email{{\tt
paola.frediani@unipv.it}}

%\date{\today}
\maketitle

\setlength{\parskip}{.1 in}

\begin{abstract}
In this paper we study the second fundamental form of the Prym map $P_{g,r}: R_{g,r} \ra {\mathcal A}^{\delta}_{g-1+r}$ in the ramified case $r>0$. We give an expression of it in terms of the second fundamental form of the Torelli map of the covering curves. 
We use this expression to give an upper bound for the dimension of a germ of a totally geodesic submanifold, and hence of a Shimura subvariety of ${\mathcal A}^{\delta}_{g-1+r}$,  contained in the Prym locus.

\end{abstract}

\section{Introduction}
Denote by $ {\mathcal R}_{g,r}$ the moduli space parametrising isomorphism classes of pairs $[(C, \alpha, R)]$ where $C$ is a smooth complex projective curve of genus $g$, $R$ is a reduced effective divisor of degree $2r$ on $C$ and $\alpha$ is a line bundle on $C$ such that $\alpha^2={\mathcal O}_{C}(R)$. To such data it is associated a double cover of $C$, $\pi:\tilde{C}\rightarrow C$ branched on $R$, with $\tilde{C}=Spec({\mathcal O}_{C}\oplus \alpha^{-1})$.

 The Prym variety
 associated to $[(C,\alpha, R)]$ is the connected component containing $0$ of the kernel of the norm map $\Nm_{\pi}: J\tilde{C} \ra JC$.
Notice that for $r>0$, $\ker \Nm_{\pi}$ is connected.  It is a  polarized abelian variety of
 dimension $g-1+r$, denoted by $P(C,\alpha,R)$ or equivalently
 $P(\tilde{C}, C)$.  The polarization $\Xi$ is induced by restricting the principal polarization on $J\tilde{C}$ and it is of type
 $
\delta=(1,\ldots, 1, \underbrace{2, ...,2}_{g \textnormal{ times }} \ )
$
for $r>0$. 
For $r=0,1$ it is twice a principal polarization and we endow $P(\tilde{C}, C)$ with this principal polarization.

 This defines the Prym map  $$P_{g,r}: {\mathcal R}_{g,r} \rightarrow {\mathcal A}^{\delta}_{g-1+r}, \ [(C,\alpha,R)] \mapsto  [(P(C,\alpha,R), \Xi)],$$
 
 where $ {\mathcal A}^{\delta}_{g-1+r}$ is the moduli space of 
 abelian varieties of dimension $g-1+r $ with a polarization of type $\delta$. We define the Prym locus as the closure in ${\mathcal A}^{\delta}_{g-1+r}$ of the image of the map $P_{g,r}$. 
 
The codifferential of $P_{g,r}$ at a generic point $[(C, \alpha, R)]$  is given by the multiplication map
\begin{equation}
\label{dp}
(dP_{g,r})^* : S^2H^0(C, K_C \otimes \alpha) \ra H^0(C, K_C^2(R))
\end{equation}

which is known to be surjective (see \cite{lo}), therefore $P_{g,r}$ is generically finite, if and only if 
$$
\dim {\mathcal R}_{g,r} \leq \dim {\mathcal A}^{\delta}_{g-1+r}.
$$
This holds if: 
either $r\geq 3$ and $g\geq 1$, or $r=2$ and $g\geq 3$, $r=1$ and $g\geq 5$,  $r=0$ and $g \geq 6$.

If $r =0$  the Prym map is
 generically injective for $g \geq 7$  (\cite{fs},
 \cite{ka}). 
If $r>0$, Marcucci and Pirola \cite{mp}, and later, for the missing cases, Marcucci and Naranjo \cite{mn} and Naranjo and Ortega \cite{no} have proved the generic injectivity in all the 
 cases except for $r=2$, $g =3$, which was previously studied  by Nagaraj and Ramanan, and also by Bardelli, Ciliberto and Verra (see \cite{nr}, \cite{bcv}) and for which 
 the degree of the Prym map is $3$.

In this paper we study the second fundamental form of the restriction of the Prym map to the open set ${\mathcal R}_{g,r}^0$ where the Prym map is an immersion, with respect to the orbifold metric  on ${\mathcal A}_{g-1+r}^{\delta}$ induced by the symmetric metric on the Siegel space ${\mathcal H}_{g-1+r}$. 

In the unramified case $r=0$, the second fundamental form of the Prym map was studied in \cite{cfprym} using the Hodge gaussian maps introduced in \cite{cpt} and the flat structure on the degree 0 line bundle $\alpha$ on  $C$. 
Here in Theorem \ref{rhotilde-rhoP} we give a description of the second fundamental form $\rho_P$ for all $r \geq 0$ in terms of the second fundamental form $\tilde{\rho}$ of the Torelli map $j: {\mathcal M}_{2g-1+r} \rightarrow {\mathcal A}_{2g-1+r}$ described in \cite{cpt}, \cite{cf2}, \cite{cfg}. 
At the point $[(C, \alpha, R)]$ we show that $\rho_P$ is obtained from $\tilde{\rho}$ by first restricting to the kernel of $(dP_{g,r})^*$ and then projecting to $S^2H^0(K_C^2(R))$. 

This also allows us to prove that the map $\rho_P$ is a lifting of the second gaussian map of the line bundle $K_C \otimes \alpha$ (see Proposition \ref{mu2}). 
This was already proved in the unramified case $r=0$ in \cite{cfprym} and previously in the case of the Torelli map in \cite{cpt}. 

In the second part of the paper we use this description of $\rho_P$ to study totally geodesic submanifolds in the Prym loci. 

We recall that a conjecture by Coleman and Oort says that for big enough genus, there should not exist Shimura subvarieties of ${\mathcal A}_g$ generically contained in the Torelli locus, i.e. contained in the Torelli locus $\overline {j({{\mathcal M}}_g)}$ and intersecting $j({\mathcal M}_g)$. 
Shimura subvarieties are totally geodesic, hence it is possible to approach the conjecture studying the second fundamental of the Torelli map. This viewpoint was used in \cite{cfg},  where an upper bound on the dimension of a germ of a totally geodesic submanifold contained in the Torelli locus, depending on $g$, is given.

In \cite{cfgp} a question similar to the one of Coleman and Oort for Pryms was asked in the case $r=0,1$, i.e. when the Prym variety is principally polarised and the Prym loci $\overline{P_{g,r}({\mathcal R}_{g,r})}$ contain the Torelli locus. 
More precisely, the question is about the existence, for big enough genus,  of Shimura subvarieties generically contained in the Prym loci.  We say that a subvariety $Z \subset {\mathcal A}_{g-1 +r}$ is generically contained in the Prym locus if $Z \subset \overline{P_{g,r}({\mathcal R}_{g,r})}$, $Z \cap P_{g,r}({\mathcal R}_{g,r}) \neq \emptyset$ and $Z$ intersects the locus of indecomposable polarised abelian varieties. 
In \cite{cfgp} examples in low dimension of Shimura curves generically contained in the Prym loci for $r=0,1$ are given, using Galois covers of $\proj^1$. 

In \cite{cf3} we gave an upper bound on the dimension of a germ of a totally geodesic submanifold contained in the Prym locus  when $r=0$, depending only on $g$, which is similar to the estimate in \cite{cfg} for the Torelli locus, that we achieved via the second fundamental form.

Here we generalise the above question for any $r \geq 0$ and we find an upper bound for the dimension of a germ of a totally geodesic submanifold contained in the Prym loci which depends on $g$ and $r$ (Theorem \ref{stima2}). This is obtained as the generic case of a bound depending also on the gonality $k$, given for a germ of a totally geodesic submanifold contained in the Prym loci passing through a point $[(C, \alpha, R)]$, with $C$ a $k$-gonal curve (Theorem \ref{stima1}).

\section{The 2nd fundamental form of the Prym map}
Let  $C$ be a smooth complex projective curve of genus $g$, $R$ a reduced divisor of degree $2r$ on $C$ and $\alpha$ a line bundle on $C$ such that $\alpha^2={\mathcal O}_{C}(R)$. To such data corresponds a double cover $\pi:\tilde{C}\rightarrow C$ branched on $R$. In the ramified case $r>0$, the Prym variety $P(C,\alpha, R)$  associated to this data is the polarised abelian variety given by the kernel of the norm map $Nm_{\pi}:J\tilde{C}\rightarrow JC$. For $r=0$ it is its connected component containing the origin. For $r>1$ the polarisation is given by $\Xi:=\Theta_{\tilde{C}|P(C,\alpha,R)}$, where $\Theta_{\tilde{C}}$ is a theta divisor of $J\tilde{C}$. For $r=0,1$ the polarisation $\Theta_{\tilde{C}|P(C,\alpha,R)}$ is twice a principal polarisation $\Xi$. Consider the Prym map $P_{g,r}: {\mathcal R}_{g,r} \rightarrow {\mathcal A}^{\delta}_{g-1+r},$ which associates to a point $[(C,\alpha,R)] \in {\mathcal R}_{g,r}$ the isomorphism class of its Prym variety $P(C,\alpha, R)$ with the polarisation $\Xi$.

As recalled in the introduction, the codifferential of the Prym map \eqref{dp} is given by the multiplication map and it is surjective at the generic point if either $r\geq 3$ and $g\geq 1$, or $r=2$ and $g\geq 3$, $r=1$ and $g\geq 5$,  $r=0$ and $g \geq 6$ (see \cite{lo}). 

In these cases we denote by ${\mathcal R}^0_{g,r}$ the non empty open subset of ${\mathcal R}_{g,r}$ where the Prym map $P_{g,r}$ is an immersion.

Consider the orbifold tangent bundle exact sequence of the Prym map
\begin{equation}
\label{tangent}
0 \rightarrow T_{{\mathcal R}^0_{g,r}} \rightarrow P_{g,r}^*T_{{\mathcal A}^{\delta}_{g-1+r}}  \rightarrow {\mathcal N_{{\mathcal R}^0_{g,r}/{\mathcal A}^{\delta}_{g-1+r}}} \rightarrow 0
\end{equation}

On ${\mathcal A}^{\delta}_{g-1+r}$ we consider the orbifold metric induced by the symmetric metric on the Siegel space ${\mathcal H}_{g-1+r}$ and the associated second fundamental form  with respect to the metric connection of the above exact sequence. 
Denote its dual by
\begin{equation}
\rho_P: {\mathcal N^*_{{\mathcal R}^0_{g,r}/{\mathcal A}^{\delta}_{g-1+r}}} \rightarrow S^2 \Omega^1_{{\mathcal R}^0_{g,r}}.
\end{equation}

To describe this second fundamental form we study the second fundamental form of the Torelli map of the covering curves $\tilde{C}$. 
Since our computations will be local, we will restrict to an open set $U$ of ${\mathcal R}^0_{g,r}$ where there is a universal family  ${\tilde{f}} :    \tilde{\mathcal {C}} \rightarrow U$ and the differential of the modular map $U \ra {\mathcal M}_{\tilde{g}}$, $[\tilde{C} \ra C] \mapsto [\tilde{C}]$ is injective.

Denote by $\tilde{H} := R^1{\tilde{f}}_*{\Co}$, $\tilde{H}_b= H^1(\tilde{C}_b, \Co)$, by ${\tilde{\mathcal F}} = \tilde{f}_*{\omega}_{\tilde{\mathcal C}|U}$ the Hodge bundle, ${\tilde{\mathcal F}}_b = H^0(\tilde{C}_b, K_{{\tilde{C}}_b})$.  
The $\Z/2\Z$ action corresponding to the $2:1$ covering gives a decomposition $\tilde{H} = \tilde{H}^+ \oplus \tilde{H}^-$ in $\pm1$ eigenspaces and an analogous decomposition $\tilde{\mathcal F} = \tilde{\mathcal F}^+ \oplus \tilde{\mathcal F}^-$. We have $\tilde{\mathcal F}^+_b \cong H^0(C_b, K_{C_b})$, $\tilde{\mathcal F}^-_b \cong H^0(C_b, K_{C_b}\otimes \alpha_b)$, $\forall b \in U$. 

The Gauss-Manin connection $\nabla_{GM}$ on $\tilde{H}$ induces a connection $\nabla^{1,0}$ on the Hodge bundle ${\tilde{\mathcal F}}$. Both connections are $\Z/2\Z$ invariant, hence we have a connection $\nabla^-$ on $\tilde{\mathcal F}^-$, and an induced connection $\nabla$ on $S^2 \tilde{\mathcal F}^-  $. 
We have the identifications on $U$:  $ \Omega^1_{{\mathcal M}_{{\tilde{g}}|U}} \cong (\tilde{f}_*{\omega}^2_{\tilde{\mathcal C}|U})$ and $\Omega^1_U $ can be identified with the subbundle preserved by the $\Z/2\Z$ action:  $\Omega^1_U = (\tilde{f}_*{\omega}^2_{\tilde{\mathcal C}|U})^+$. Note in fact that at a point $b:=[\tilde{C} \ra C] \in U$, ${\Omega^1_U}_b = H^0(C, K_C^2(R)) \cong H^0(K_{\tilde{C}}^2)^+$, as one can easily check via the projection formula.
Moreover  $P^*\Omega^1_{{{\mathcal A}^{\delta}_{g-1+r}}} \cong S^2\tilde{\mathcal F}^-$ and the connection $\nabla$ corresponds to the connection associated to the Siegel metric on ${{{\mathcal A}^{\delta}_{g-1+r}}}$.

Consider the exact sequence
\begin{equation}
\label{II_U}
0 \ra {\mathcal I}_2 \ra S^2{\tilde{\mathcal F}}  \stackrel {m}\rightarrow  \tilde{f}_*{\omega}^2_{\tilde{\mathcal C}|U}\ra 0
\end{equation}
where the map $m$ is the multiplication map and it is the dual of the differential of the Torelli map of the curves $\tilde{C}$ on $U$.

The second fundamental form of the exact sequence {\eqref{II_U}} is a map 

\begin{equation} 
\Psi: {\mathcal I}_2 \ra   \tilde{f}_*{\omega}^2_{\tilde{\mathcal C}|U} \otimes \Omega^1_U \cong \tilde{f}_*{\omega}^2_{\tilde{\mathcal C}|U} \otimes  (\tilde{f}_*{\omega}^2_{\tilde{\mathcal C}|U})^+
\end{equation}

Since the multiplication map $m$ is $\Z/2\Z$ equivariant, we also have the exact sequence:

\begin{equation} 
0 \ra {\mathcal I}_2^+ \ra (S^2{\tilde{\mathcal F}})^+  \stackrel {m} \rightarrow  (\tilde{f}_*{\omega}^2_{\tilde{\mathcal C}|U})^+\ra 0
\end{equation} 

Clearly we have $(S^2{\tilde{\mathcal F}})^+ \cong S^2{\tilde{\mathcal F}}^+ \oplus S^2{\tilde{\mathcal F}}^-$ and the restriction of the multiplication map $m$ to $S^2{\tilde{\mathcal F}}^-$ is the dual of the differential of the Prym map $P$ (see \cite{lo}). 
More precisely, with the identifications  on $U$: 
 $P^*\Omega^1_{{{\mathcal A}^{\delta}_{g-1+r}}} \cong S^2\tilde{\mathcal F}^-, \ \Omega^1_U = (\tilde{f}_*{\omega}^2_{\tilde{\mathcal C}|U})^+$, the dual of the exact sequence \eqref{tangent} on $U$ can be written as

\begin{equation}
\label{G}
0 \ra  {\mathcal G} \ra  S^2{\tilde{\mathcal F}}^- \stackrel {m} \ra  (\tilde{f}_*{\omega}^2_{\tilde{\mathcal C}|U})^+ \ra 0
\end{equation}
where ${\mathcal G} = S^2{\tilde{\mathcal F}}^-  \cap {\mathcal I}_2^+ {{{\mathcal A}^{\delta}_{g-1+r}}} \cong {\mathcal N}^*_{U/{{{\mathcal A}^{\delta}_{g-1+r}}}}$. 
So the dual of the second fundamental form of the Prym map is a map 
\begin{equation}
\rho_P: {\mathcal G} \ra  (\tilde{f}_*{\omega}^2_{\tilde{\mathcal C}|U})^+ \otimes  (\tilde{f}_*{\omega}^2_{\tilde{\mathcal C}|U})^+,
\end{equation} 
which is symmetric, since it is the second fundamental form of an immersion.
Clearly we have 
\begin{equation}
\label{ro_P}
\rho_P = p \circ \Psi_{| {\mathcal G}}
\end{equation}
where $p:  (\tilde{f}_*{\omega}^2_{\tilde{\mathcal C}|U})\otimes   (\tilde{f}_*{\omega}^2_{\tilde{\mathcal C}|U})^+ \ra S^2 (\tilde{f}_*{\omega}^2_{\tilde{\mathcal C}|U})^+$ is the natural projection.

Denote by $\tilde{g}:= 2g-1+r$ and consider now the Torelli map $\tilde{j}: {\mathcal M}^0_{\tilde{g}} \ra {\mathcal A}_{\tilde{g}}$, where $ {\mathcal M}^0_{\tilde{g}}$ is the complement of the hyperelliptic locus. Then $\tilde{j}$ is an immersion and we denote by $\tilde{\rho}$ the dual of the second fundamental form of $\tilde{j}$. Since we are working in a local setting,  we can assume that we have a modular map $\mu: U \ra V $ where $V$ is  an open subset of ${\mathcal M}^0_{\tilde{g}}$ on which  there exists a universal family $\phi: \tilde{\mathcal C}_V \ra V$. The family $\tilde{\mathcal C}$ is the pullback of $\tilde{\mathcal C}_V$ via the map $\mu$. 

On $V$ we have 
\begin{equation}
0 \ra  {\mathcal I}_2(\omega_{\tilde{\mathcal C}_V|V}) \ra S^2(\phi_* \omega_{\tilde{\mathcal C}_V|V}) \stackrel {m} \ra  \phi_* \omega^2_{\tilde{\mathcal C}_V|V} \ra 0
\end{equation} 
and the multiplication map $m$ is the dual of the differential of the Torelli map. 

On $V$ the dual of the second fundamental form of the Torelli map is a map
\begin{equation} 
\tilde{\rho}:  {\mathcal I}_2(\omega_{\tilde{\mathcal C}_V|V}) \ra \phi_* \omega^2_{\tilde{\mathcal C}_V|V} \otimes \phi_* \omega^2_{\tilde{\mathcal C}_V|V} 
\end{equation}

The pullback of the above exact sequence on $U$ via $\mu$ is the exact sequence {\eqref{II_U}}. 
Hence we have 
\begin{equation}
\label{PSI}
\Psi =  q \circ \mu^* \tilde{\rho}
\end{equation}
where $q:   (\tilde{f}_*{\omega}^2_{\tilde{\mathcal C}|U})\otimes   (\tilde{f}_*{\omega}^2_{\tilde{\mathcal C}|U}) \ra   (\tilde{f}_*{\omega}^2_{\tilde{\mathcal C}|U})\otimes   (\tilde{f}_*{\omega}^2_{\tilde{\mathcal C}|U})^+ $ is the natural projection.

We have the following

\begin{TEO}
\label{rhotilde-rhoP}
The dual of the second fundamental form of the Prym map on $U$ is obtained as 
$\rho_P=  p' \circ (\mu^* \tilde{\rho})_{| {\mathcal G}} $, where $p' : S^2(\tilde{f}_*\omega^2_{{\mathcal C}|U}) \ra S^2((\tilde{f}_*\omega^2_{{\mathcal C}|U})^+)$ is the natural projection.
\end{TEO} 
\proof
By \eqref{ro_P} and \eqref{PSI} we have $\rho_P=  p \circ \Psi_{| {\mathcal G}}=  p \circ q \circ (\mu^* \tilde{\rho})_{| {\mathcal G}} = p' \circ  (\mu^* \tilde{\rho})_{| {\mathcal G}}$. 
\qed

\

At the point $b_0:= [(C,\alpha,R)] \in U$ corresponding to the $2:1$ cover $\pi:\tilde{C} \ra C$, the space $P_{g,r}^*\Omega^1_{{{\mathcal A}^{\delta}_{g-1+r}}, b_0}$ is isomorphic to $S^2H^0(K_C \otimes \alpha)$,    $\Omega^1_{{\mathcal R}^0_{g,r},b_0}$ is isomorphic to $ H^{0}(K_C^2(R))$, ${\mathcal G}_{b_0} \cong I_2(K_C \otimes \alpha)$, and the dual of the exact sequence \eqref{tangent} at the point $b_0$, that is the exact sequence \eqref{G} at $b_0$ becomes

$$0 \rightarrow I_2(K_C\otimes \alpha)  \rightarrow S^2H^0(K_C \otimes \alpha) \stackrel {m}\rightarrow H^0(K_C^2(R))
\rightarrow 0. $$

The dual of the second fundamental form of the Prym map at the point $b_0$  is a map

 \begin{equation}
 \label{II}
 \rho_P: I_2(K_C\otimes \alpha)\rightarrow S^2H^0(K_C^2(R))
\end{equation}

Observe that $\forall Q \in I_2(K_C\otimes \alpha)  \stackrel{\pi^*} \hookrightarrow  I_2(K_{\tilde{C}})^+$, $\forall v_1, v_2 \in H^1(T_C(-R)) \cong H^1(T_{\tilde{C}})^+$, by Theorem \eqref{rhotilde-rhoP} we have: 
\begin{equation}
\label{rho(Q)}
\rho_P(Q)(v_1  \odot v_2) = \tilde{\rho}(\pi^*Q)(v_1 \odot v_2).
\end{equation} 

Denote by 
\begin{equation}
\label{mutilde}
\tilde{\mu} _2: I_2(K_{\tilde{C}}) \ra H^0(K^4_{\tilde{C}})
\end{equation}
the second Gaussian map of the canonical bundle $K_{\tilde{C}}$ and by 

\begin{equation}
\label{mu2}
{\mu}_2:= \mu_{2,K_C \otimes \alpha}  : I_2(K_{{C}} \otimes \alpha) \ra H^0(K^4_{{C}}(R))
\end{equation}
the second Gaussian map of the line bundle $K_C \otimes \alpha$ (for the definition and a local expression of the second Gaussian maps see  e.g. \cite{cf1}, section 2). 
Notice that $\tilde{\mu} _2$ is equivariant, hence it induces a map $\tilde{\mu} _2:I_2(K_{\tilde{C}})^+ \ra H^0(K^4_{\tilde{C}})^+ \cong H^0(K^4_C(2R))$. 

We have the following  

\begin{LEM}
\label{mu}
For every $Q \in I_2(K_C \otimes \alpha) $, $\mu_2(Q)  = \tilde{\mu} _2(\pi^*Q) $ via the inclusion 
$H^0(K^4_{{C}}(R)) \subset H^0(K^4_C(2R))$.

\end{LEM}
\proof

We show the equality by a local computation. 
For a point $P \not \in R$, take local coordinates $z$ in a neighbourhood $V$ of  $P$ and $w$ in a neighbourhood $U$ of a point $T \in \pi^{-1}(P)$ such that the local expression of $\pi:U \ra V$ is $w \mapsto w=z$. Since $\alpha^2 = \OO_C(R)$, $\alpha^2_{|V} = \OO_V$ and we choose a local frame $a$ of $\alpha$ on $V$ such that $a^2 =1$ and  $(\pi^*a)_{|U} = 1$. 
Fix  a basis $\{\omega_i\}$ of $H^0(K_C \otimes \alpha)$, so locally $\omega_i = f_i(z) dz \otimes a$. Then on $U$ we have $\pi^*(\omega_i) = f_i(w)dw$. 
 Take a quadric $Q =  \sum_{i,j}  a_{i,j} \omega_i \odot \omega_j \in I_2(K_C \otimes \alpha) $, hence locally  $\pi^*Q = \sum_{i,j}  a_{i,j} f_i(w)dw \odot f_j(w)dw$ and we have 
$\tilde{\mu} _2 (\pi^*Q) = -\sum_{i,j}  a_{i,j} f'_i(w)f'_j(w)(dw)^4. $
On the other hand on $V$, $\mu_2(Q) = -\sum_{i,j} a_{i,j} f'_i(z) f'_j(z) (dz)^4 a^2  $, hence the statement follows, since $a^2=1$. 

Observe that the equality can also be checked locally around a critical point $T$ over a point $P \in R$. So we can assume that the map $\pi:U \ra V$ is  of the form $w \mapsto w^2=z$. 
Now $\alpha^2_{|V} = \OO_V(P)$ and  $(\pi^*\alpha)_{|U} = \OO_U(T)$ so we choose a local frame $a$ of $\alpha$ on $V$ such that $a^2 =\frac{1}{z}$ and  $(\pi^*a)_{|U} = \frac{1}{w}$. 
Now locally $\pi^*(\omega_i) = 2f_i(w^2)dw$ and $\pi^*(Q) = 4\sum_{i,j}  a_{i,j} f_i(w^2)dw \odot f_j(w^2)dw$. So we have 
$$\tilde{\mu} _2 (\pi^*Q) = -4\sum_{i,j}  a_{i,j} f'_i(w^2)f'_j(w^2)4w^2(dw)^4= -\sum_{i,j}  a_{i,j} f'_i(z)f'_j(z)\frac{(dz)^4}{z}. $$
On the other hand locally  $\mu_2(Q) = -\sum_{i,j} a_{i,j} f'_i(z) f'_j(z) (dz)^4 a^2 =  \tilde{\mu} _2 (\pi^*Q) $. 
\qed

We recall the definition of a Schiffer variation of a line bundle $L$ on a curve $C$  at a point $p  \in C$.
Consider the evaluation map $v: H^0(K_C \otimes L) \otimes \OO_C \ra K_C\otimes L $  and its dual $$\xi_L: T_C \otimes L^{-1} \ra H^1(L^{-1}) \otimes \OO_C.$$
A Schiffer variation in $P$ is a generator of the image of the map $$(T_C \otimes L^{-1})_P \ra H^1(L^{-1}).$$
Note that this map can be seen as the coboudary map of the exact sequence 
\begin{equation}
\label{schiffer}
0 \ra L^{-1} \ra L^{-1}(P) \ra  L^{-1}(P)_{|P} \ra 0,
\end{equation}
with the identification ${T_C}_{|P}= \OO_C(P)_{|P}$.

If we fix a local coordinate $z$ centred in $P$, we choose the Schiffer variation  $\xi_{P,L}$ in $P$ to be the Doulbeault class of $\theta_{P,L} = \frac{\bar{\partial}b_P}{z} \otimes l^{-1} \in A^{0,1}(L^{-1})$, where $b_P$ is a bump function at $P$ which is equal to one in a neighborhood of $P$ and $l$ is a local frame of $L$.

\begin{PROP}
\label{mu2}
We have the following commutative diagram

$$\xymatrix{
 I_2(K_C\otimes \alpha) \ar[d]^{-4\pi i \mu_2}\ar[r]^{\rho_P} & S^2(H^0(K^2_C(R))) \ar[d]^{m} \\
H^0(K^4_C(R)) \ar[r]& H^0(K^4_C(2R))}
$$

\end{PROP}
\proof
We recall that we have a similar statement for $\tilde{\rho}$ and $\tilde{\mu_2}$, namely $m \circ \tilde{\rho} = -2 \pi i  \tilde{\mu}_2$  (see {\cite[Thm. 3.1]{cpt}, \cite[Thm. 2.2]{cfg}}).
Denote as usual by $\pi: \tilde{C} \ra C$ the double cover and by $\sigma$ the covering involution on $\tilde{C}$, take a point $T \in \tilde{C}$ with $T \neq \sigma(T)$ and fix a local coordinate on $\tilde{C}$ around $T$ (and correspondingly around $\sigma(T)$)  and on $C$ around $P:= \pi(T)$.  For a point $S \in \tilde{C}$, denote by $\xi_S:= \xi_{S, K_{\tilde{C}}}$.

Then by  {\cite[Thm. 3.1]{cpt}, \cite[Thm. 2.2]{cfg}} we have: 
\begin{equation}
\label{annullamento}
\tilde{\rho}(\pi^*Q)(\xi_T \odot \xi_{\sigma(T)}) = -4 \pi i (\pi^*Q)(T, \sigma(T)) \cdot \tilde{\eta}_T(\sigma(T)) = -4 \pi i Q(P,P) \cdot \tilde{\eta}_T(\sigma(T)) =0,
\end{equation}
since $Q(P,P)=0$ and $\tilde{\eta}_T \in H^0(K_{\tilde{C}}(2T))$ has only one double pole in $T$, so it is  holomorphic around $\sigma(T)$. 
Set $v:= \xi_T + \xi_{\sigma(T)} \in H^1(T_{\tilde{C}})^+$. By \eqref{annullamento} and by the $\Z/2\Z$ equivariance of $\tilde{\rho}$ we have: 
\begin{equation}
\label{1}
2\tilde{\rho}(\pi^*Q)(\xi_T \odot \xi_T)= \tilde{\rho}(\pi^*Q)(v\odot v) = \rho_P(Q)(v \odot v)
\end{equation}
where the last equality is \eqref{rho(Q)}. 

By  {\cite[Thm. 3.1]{cpt}, \cite[Thm. 2.2]{cfg}} we have 
\begin{equation}
\label{2}
\tilde{\rho}(\pi^*Q)(\xi_T \odot \xi_T) = m(\tilde{\rho}(\pi^*(Q))(T) = -2 \pi i \tilde{\mu}_2(\pi^*Q)(T) = -2 \pi i \mu_2(Q)(P),
\end{equation}
by Lemma \eqref{mu}.

Consider the line bundle $L = K_C(R)$, fix a local coordinate $z$ in $P$ and the corresponding  Schiffer variation $\xi_{P,L} \in H^1(T_C(-R))$.  Choose in $T$ and $\sigma(T)$ the  local coordinates determined by $z$ and take  the associated Schiffer variations $\xi_T$, $\xi_{\sigma(T)}$ of $K_{\tilde{C}}$ at $T$ and $\sigma(T)$. Then with the identification of $H^1(T_C(-R))$ with $H^1(T_{\tilde{C}})^+$, the Schiffer variation $\xi_{P,L}$ corresponds to $v =\xi_T + \xi_{\sigma(T)} $. 
This can be checked as follows:  take the exact sequence 
\begin{equation}
\label{T}
0 \ra T_{\tilde{C}} \ra T_{\tilde{C}}(T + \sigma(T)) \ra T_{\tilde{C}}(T + \sigma(T))_{|T + \sigma(T)} \ra 0.
\end{equation}
The image of its coboundary map $H^0(T_{\tilde{C}}(T + \sigma(T))_{|T + \sigma(T)}) \ra H^1(T_{\tilde{C}})$ is the subspace generated by $\{\xi_T, \xi_{\sigma(T)}\}$ and the $\Z/2\Z$-invariant part of this image is generated by $v=\xi_T + \xi_{\sigma(T)}$. 
Since $\pi$ is finite,  $R^1\pi_*  T_{\tilde{C}}  =0$,  so, if we apply $\pi_*$ to \eqref{T}, using the projection formula and the isomorphism $\pi_*\OO_{\tilde{C}} \cong \OO_C \oplus \alpha^{-1}$, we get the exact sequence
\begin{equation}
0 \ra (T_C \otimes \alpha^{-1} ) \oplus  T_{C}(-R) \ra  (T_C \otimes \alpha^{-1}(P) )\oplus T_{C}(-R)(P) \ra  \pi_* ( T_{\tilde{C}}(T + \sigma(T))_{|T + \sigma(T)}) \ra 0.
\end{equation}

Taking the invariant part, we get the exact sequence \eqref{schiffer} for $L = K_C(R)$:
\begin{equation}
\label{P}
0 \ra T_{C}(-R) \ra T_{C}(-R)(P) \ra  T_{C}(-R)(P)_{|P} \ra 0.
\end{equation}
In particular we have the isomorphism $H^0((\pi_*(T_{\tilde{C}}(T + \sigma(T))_{|T + \sigma(T)}))^+ )\cong H^0( T_{C}(-R)(P)_{|P} )$. 
Recall that the Schiffer variation $\xi_{P,L}$ is a generator of the image of the coboundary map $H^0( T_{C}(-R)(P)_{|P} ) \ra H^1( T_{C}(-R) ) \cong H^1(T_{\tilde{C}})^+$ of the exact sequence \eqref{P}. 

%The coboundary map of \eqref{T} can be identified with the coboundary map
%$H^0(\pi_*(T_{\tilde{C}}(T + \sigma(T))_{|T + \sigma(T)} ))  \cong H^0(T_{\tilde{C}}(T + \sigma(T))_{|T + \sigma(T)}) \ra H^1(\pi_*( T_{\tilde{C}} ) )\cong H^1(T_{\tilde{C}})$ of the exact sequence obtained applying $\pi_*$ to the exact sequence \eqref{T}.

%Since $\xi_T + \xi_{\sigma(T)}$ is invariant by the $\Z/2\Z$ action, it  is a generator of the  one dimensional subspace of image of the restriction of the above coboundary map  to the subspace 
%$H^0((\pi_*(T_{\tilde{C}}(T + \sigma(T))_{|T + \sigma(T)}))^+ )\cong H^0( T_{C}(-R)(P)_{|P} )$. 
Form this it follows that the subspace $\langle v  \rangle \subset H^1(T_{\tilde{C}})^+  \cong H^1(T_C(-R))$ is identified with the subspace generated by $\xi_{P,L}$. 
By the choice we made of the coordinates around $P$, $T$ and $\sigma(T)$ one can check the identification of  $v$ with $\xi_{P,L}$.

So we have  
\begin{equation}
\label{3}
\rho_P(Q)(v \odot v) = \rho_P(Q)(\xi_{P,L} \odot \xi_{P,L}) = (m(\rho_P(Q)))(P),
\end{equation}
by the definition of $\xi_{P,L}$, 
and putting together \eqref{1}, \eqref{2}, \eqref{3} we get  
\begin{equation}
(m(\rho_P(Q)))(P) = -4 \pi i \mu_2(Q)(P) 
\end{equation}
for all $P$ which is not a critical value of $\pi$. Hence $m \circ \rho_P(Q)$ and $-4 \pi i \mu_2(Q)$ are two sections of $H^0(K^4_C(R))$ that coincide in the complement of a finite set in $C$, so they coincide.  
\qed

 \section{Totally geodesic submanifolds}

In this section, following the ideas of \cite{cfg} and \cite{cf3}, we give an upper bound for the dimension of a germ of a totally geodesic submanifold of $ {\mathcal A}^{\delta}_{g-1+r}$ contained in $P_{g,r}( {\mathcal R}^0_{g,r})$.

\begin{PROP}
  \label{rank}
  Assume that $[(C,\alpha, R)] \in {\mathcal R_{g,r}^0}$ is such that $C$ is a  $k$-gonal curve of genus $g$, with $g +r\geq k +3$
   and $\alpha^2={\mathcal O}_C(R)$.
  \begin{enumerate}
  \item If $r > k+1$,  then there exists a quadric $Q \in I_2(K_C \otimes \alpha)$ such that $\operatorname{rank}\rho(Q) \geq  2g-2 -k+r$.
  \item If $r \leq k+1$,  then there exists a quadric $Q \in I_2(K_C \otimes \alpha)$ such that $\operatorname{rank}\rho(Q) \geq  2g-2k -4+2r$.

  \end{enumerate}
 \end{PROP}
\proof Let $F$ be a line bundle on $C$ such that $|F|$ is
a $g^1_k$ and choose a basis $\{x,y\}$ of $H^0(F)$.  Set $M = K_C \otimes \alpha \otimes F^{-1}$ and denote by $B$ the base locus of $|M|$. By Riemann Roch
\begin{equation}
\label{RR}
h^0(M)= h^0(M(-B))=h^0(F\otimes \alpha^{-1})+g-1-k+r\geq g-k+r-1\geq 2\quad
\end{equation}
by assumption.

Note that in case (1) $B = \emptyset$, since $deg(M)=2g-2+r-k>2g-1$.

Take a pencil $\langle t_1,t_{2}\rangle$ in $H^0(M)$.  If $B \neq \emptyset$, write $t_i = t'_i s$ for a
section $s\in H^0(C,\OO_C(B)) $ with $\operatorname{div} (s) = B$.
Then $\langle t'_1, t'_2\rangle $ is a base point free pencil in $|M(-B)|$.  Let
$\psi : C \rightarrow  \proj^1$ be the morphism induced by this pencil and $\tilde{\psi}=\psi\circ \pi: \tilde{C}\rightarrow \proj^1$ and set $d:= deg(\psi)= deg(M(-B))$. 
Denote by $\phi$ the morphism induced
by the pencil $|F|$ and  and $\tilde{\phi}=\phi\circ \pi: \tilde{C}\rightarrow \proj^1$.

Consider the rank 4 quadric $Q:= xt_1 \odot yt_2 - xt_2 \odot
yt_1$. Clearly $Q \in I_2(K_C \otimes \alpha)$. We want to show that $rk\tilde{\rho}(\pi^*Q)\geq d$.

% As in the proof of \cite{cfg} Theorem 4.1 one can  show that if
 % $\{P_1,...,P_d\}$ is a fibre of the morphism $\psi$ over a regular
%value.
%then the Schiffer variations $\xi_{P_1},...,
%\xi_{P_d}$ are linearly independent in $H^1(C,T_C)$.
Consider the set $E:= \psi(R\cup Crit(\phi) \cup
Crit(\psi) \cup B)$ where $Crit(\phi) $ (resp.  $Crit(\psi) $)
denote the set of critical points of $\phi$ (resp. $\psi$).  Let $z
\in {\proj}^1 \setminus E$ and let
$\{P_1,\ldots,P_d\}$ be the fibre of
$\psi$ over $z$. By changing coordinates on $\proj^1$ we can assume $z =
[0,1]$, i.e. $t'_1(P_i) = 0$ for $i =1,\ldots d$. Then clearly
$t_1(P_i) = 0$, so $Q(P_i,P_j) = 0$ for all $i,j$.
Set $\{T_i,\sigma(T_i)\}=\pi^{-1}(P_i)$, so $\pi^*Q(T_i,T_j)=\pi^*Q(T_i,\sigma(T_j))=Q(P_i,P_j) =0.$

Let us fix a local coordinate at the relevant points and  write $\xi_T:= \xi_{T, K_{\tilde C}}$ for a Schiffer
variation of $\tilde{C}$ at $T$. 
Set $v_i:=\xi_{T_i}+\xi_{\sigma(T_i)}$. Clearly $v_i \in H^1(T_{\tilde C})^+\simeq H^1(T_{C}(-R))$, so by \eqref{rho(Q)} we have ${\rho}_P(Q)(v_i\odot v_j)=\tilde{\rho}(\pi^*Q)(v_i\odot v_j)$.

Hence, by  \cite[Thm. 2.2]{cfg}, for $i\neq j$
\begin{equation}
\label{MM}
\begin{gathered}
{\rho}_P(Q)(v_i\odot v_j)=\tilde{\rho}(\pi^*Q)(v_i\odot v_j)=2\tilde{\rho}(\pi^*Q)(\xi_{T_i}\odot \xi_{T_j})+2\tilde{\rho}(\pi^*Q)(\xi_{\sigma(T_i)}\odot \xi_{T_j})=\\
-8 \pi i (\pi^*Q)(T_i,T_j)\tilde{\eta}_{T_j}(T_i)-8 \pi i (\pi^*Q)(\sigma(T_i),T_j)\tilde{\eta}_{T_j}(\sigma(T_i))=0
\end{gathered}
\end{equation}
and 
\begin{equation}
\label{MM1}
\begin{gathered}
\tilde{\rho}(\pi^*Q)(v_i\odot v_i)=2\tilde{\rho}(\pi^*Q)(\xi_{T_i}\odot \xi_{T_i})+2\tilde{\rho}(\pi^*Q)(\xi_{\sigma(T_i)}\odot \xi_{T_i})=\\
-4 \pi i \tilde{\mu_2}(\pi^*Q)(T_i)-8 \pi i (\pi^*Q)(\sigma(T_i),T_i)\tilde{\eta}_{T_i}(\sigma(T_i))=-4 \pi i \tilde{\mu_2}(\pi^*Q)(T_i) = -4 \pi i \mu_2(Q)(P_i).
\end{gathered}
\end{equation}

 For a rank 4 quadric the second
Gaussian map can be computed as follows: $ \mu_2(Q) = \mu_{1,{F}}(x
\wedge y)\mu_{1,M}(t_1 \wedge t_2)$, where $\mu_{1,F}$ and $\mu_{1,M}$ are the first Gaussian maps of the line bundles $F$ and $M$ (see \cite[Lemma 2.2]{cf1}).  Now
$\mu_{1,F}(x\wedge y) (P_i) \neq 0$, because $P_i \not \in
Crit({\phi})$ by the choice of $z$.  Moreover
$P_i $ is not in the base locus $B$.  On $C\setminus B$ the morphism $\psi$ coincides
with the map associated to $\langle t_1, t_2\rangle$. Since $P_i \not \in
Crit(\psi)$, it is not a critical point for the latter map. Therefore
also $\mu_{1,M}(t_1 \wedge t_2) (P_i) \neq 0$.
Thus $\mu_2(Q)(P_i) = \mu_{1,F}(x \wedge y)(P_i)\mu_{1,M}(t_1 \wedge
t_2)(P_i) \neq 0$ for every $i=1,..., d$.

We claim that the vectors $\{v_1,...,v_d\}$ are linearly independent in $H^1(T_{\tilde C})^+\simeq H^1(T_{C}(-R))$.

In fact we show that the subspace  $W:= \langle \xi_{T_1}, \xi_{\sigma(T_1)}, ...,\xi_{T_d}, \xi_{\sigma(T_d)} \rangle$ of  $H^1(\tC, T_{\tC})$ has dimension $2d$. This is equivalent to say that the annihilator $Ann(W)$ of $W$ in $H^0(\tC, K^2_{\tC})$ has codimension $2d$. Observe that $Ann(W) = H^0(\tC, K^2_{\tC}(-D))$, where $D = T_1 + \sigma(T_1) + ...+ T_d + \sigma(T_d)$. Then by Riemann Roch we have $h^0(K^2_{\tC}(-D)) = h^0(T_{\tC}(D)) + 4(\tilde{g}-1) - 2d -(\tilde{g}-1) = 3\tilde{g}-3 -2d$, since $deg(T_{\tC}(D)) = -2(\tilde{g}-1) + 2d = -4g -2r + 4 + 2d \leq -4g -2r +4 + 2(2g-2 + r -k-deg(B)) =  -2k -2deg(B)<0$.

This shows that $\{ \xi_{T_1}, \xi_{\sigma(T_1)}, ...,\xi_{T_d}, \xi_{\sigma(T_d)} \}$ are linearly independent  in $H^1(\tC, T_{\tC})$ and hence also $\{v_1,...,v_d\}$ are linearly independent in $H^1(T_{\tilde C})^+$.

 By \eqref{MM},\eqref{MM1}  one immediately obtains that the restriction of
$\tilde{\rho}(Q)$ to the subspace $ W':= \langle v_1,...,v_d \rangle$ is represented in the basis
$\{v_1,...,v_d\}$ by a diagonal matrix with entries 

$-4 \pi i \mu_2(\pi^*Q)(T_i) = -4 \pi i \mu_2(Q)(P_i)\neq 0$ on the diagonal.   So
$\rho_P(Q)$ has rank at least $d$.

 In case $(1)$ $B$ is empty, hence $d = 2g-2 +r -k$.  In case $(2)$, by Clifford Theorem we have:
 $$2(h^0(M(-B)) - 1) \leq deg(M(-B)) = 2g-2 +r- k - deg(B),$$
 hence
 $$deg(B) \leq 2g-2+r -k - 2h^0(M(-B)) +2 \leq 2g  +r -k - 2(g-k+r-1) = k -r +2 $$
 and
 $d= 2g-2+r-k-
\deg(B)\geq 2g+2r-2k-4$.

\qed

\begin{TEO}
  \label{stima1}
  Assume that $[(C,\alpha,R)] \in {\mathcal R}^0_{g,r}$ where $C$ is a $k$-gonal curve of genus $g$ with $g+r\geq k+3$. Let $Y$ be a germ of a totally geodesic submanifold of
  ${\mathcal A}^{\delta}_{g-1+r}$ which is contained in $P_{g,r}( {\mathcal R}^0_{g,r})$ and passes through
  $P(C,\alpha,R)$. Then 
   \begin{enumerate}
  \item If $r > k+1$,  then $\dim Y\leq 2g-2+\frac{3r+k}{2}$.
  \item If $r \leq k+1$, then $\dim Y\leq 2g+r+k-1$.

  \end{enumerate}
\end{TEO}
\proof Since $Y$ is totally geodesic, for any $v \in T_{[(C,\alpha,R)]}Y$ we must have that $\rho(Q)(v \odot v) = 0$ for any $Q$ in $I_2(K_C\otimes\alpha)$. Hence if a quadric $Q$ is such that the rank of $\rho(Q)$ is at least $m$
 \begin{gather*}
  \dim T_{[C]} Y \leq (3g-3 +2r) - \frac{m }{2}
  .
\end{gather*}
The result then follows by the existence of a quadric $Q
\in I_2(K_C\otimes \alpha )$ shown in Theorem \ref{rank}, with  $rank (\rho(Q))\geq m$, where $m=2g-2 -k+r$ in case (1) and $m=2g-2k-4+2r$ in case (2).
\qed
\begin{remark}
  In case (2) of Proposition \ref{rank} if $|M|$ is base point free we have the same estimate as in case (1), namely $ \rank \rho(Q) \geq 2g-2-k+r$.  So in this case the bound on the dimension of  a
  germ of a totally geodesic submanifold  $Y$  contained in the
  Prym locus and passing through $P(C, \alpha,R)$ with $C$ a $k$-gonal curve such that $|M|$ is base point free becomes:
  $\dim Y \leq 2g-2+\frac{3r+k}{2}$.
\end{remark}

We will now give a bound on the dimension of a germ a totally geodesic submanifold of
  ${\mathcal A}^{\delta}_{g-1+r}$ contained in $P_{g,r}( {\mathcal R}^0_{g,r})$ which does not depend on the gonality.

\begin{TEO}
  \label{stima2}
  Let $Y$ be a germ of a totally geodesic submanifold of
  ${\mathcal A}^{\delta}_{g-1+r}$ which is contained in $P_{g,r}( {\mathcal R}^0_{g,r})$.
    \begin{enumerate}
  \item If $g<2r-5$,  then $\dim Y\leq \frac{9}{4}g+\frac{3}{2}r-\frac{5}{4}$.
  \item If $g \geq 2r-5$, then $\dim Y\leq  \frac{5}{2}g+r+\frac{1}{2}$.

  \end{enumerate}
  
\end{TEO}
\proof 
Recall that the gonality $k$ of a genus $g$ curve is at most $ [(g+3)/{2}]$. 

If $g <2r-5$, then $k \leq [(g+3)/{2}]< r-1$,  hence  statement $(1)$ follows immediately by Theorem \ref{stima1} $(1)$. 

If $g \geq 2r-5$, $r \leq [(g+3)/{2}] +1$.
 If $r \leq k+1$, we apply inequality $(2)$ of Theorem \ref{stima1} and we get statement $(2)$. 
If $k <[(g+3)/{2}]$ and $r >k+1$, then estimate $(1)$ of Theorem \ref{stima1} applies, so $\dim Y \leq 2g-2+\frac{3r+k}{2} < \frac{9}{4}g + \frac{3}{2} r - \frac{5}{4} \leq \frac{5}{2}g+r+\frac{1}{2}$.
 \qed

\begin{remark}
In \cite{cfgp} examples of Shimura curves (hence totally geodesic) of ${\mathcal A}_{g-1 +r}$  contained in $P_{g,r}({\mathcal R}_{g,r})$ when $r =0,1$ have been constructed using families of Galois covers of $\PP^1$. The examples when $r=0$ are all contained in ${\mathcal A}_{g-1}$ with $g \leq 13$, while the ones with $r=1$ are all contained in ${\mathcal A}_{g}$ with $g \leq 8$. 
\end{remark}


\begin{thebibliography} {99}
\bibitem{bcv}
 Bardelli,  ~F., Ciliberto,  ~C., Verra,  ~A., Curves of minimal genus on a general abelian variety,  Compos. Math. 96 (1995), no. 2, 115--147. 


%\bibitem{b} A. Beauville, {\em{Prym varieties and the Schottky
%problem}}, Inventiones Math. \textbf{41} (1977), 149-96.



  \bibitem{cf1}Colombo, ~E., Frediani, ~P.,  Some results on the second Gaussian map for curves. Michigan Math.
 J. Vol. 58,  3 (2009), 745-758.
 
 \bibitem{cf2}Colombo, ~E., Frediani, ~P., Siegel metric and curvature of the moduli space of curves.
 Transactions of the Amer. Math. Soc. 362 (2010), no. 3, 1231-1246.

\bibitem{cf3}Colombo, ~E., Frediani  ~P.,  A bound on the dimension of a totally geodesic submanifold in the Prym locus. Collectanea Mathematica. DOI: 10.1007/s13348-018-0215-0.  

\bibitem{cfprym}Colombo, ~E., Frediani, ~P., Prym map and second Gaussian map for Prym-canonical line bundles. Adv. Math. 239 (2013), 47-71.

\bibitem{cfg} Colombo, ~E., Frediani, ~P., Ghigi, ~A., On totally geodesic submanifolds in the Jacobian locus. Internat. J. Math. 26 (2015), no. 1, 1550005, 21 pp.


\bibitem{cfgp} Colombo, ~E., Frediani, ~P., Ghigi, ~A., Penegini ~M., Shimura curves in the Prym locus.
Communications in Contemporary Mathematics. DOI: 10.1142/S0219199718500098. 

\bibitem{cpt} Colombo,~E., Pirola,~G.P., Tortora,~A., Hodge-Gaussian
maps, Ann. Scuola Normale Sup. Pisa Cl. Sci. (4) {\bf 30} (2001),
no. 1, 125-146.
%\bibitem{fl} Farkas, Gavril; Ludwig, Katharina, The Kodaira dimension of the moduli space of
%Prym varieties. J. Eur. Math. Soc. (JEMS) 12 (2010), no. 3,
%755-795.
\bibitem{fs}  Friedman, Robert; Smith, Roy,  The generic Torelli theorem for the Prym map. Invent. Math. 67 (1982), no. 3, 473Ð490.
%\bibitem{green}Green,~M.~L., Infinitesimal methods in Hodge theory, in {\em Algebraic Cycles and Hodge Theory}, Torino 1993,
%Lecture Notes in Mathematics, 1594. Springer, Berlin, (1994),
%1-92.
\bibitem{ka} Kanev, V. I., A global Torelli theorem for Prym varieties at a general point. (Russian)
Izv. Akad. Nauk SSSR Ser. Mat. 46 (1982), no. 2, 244-268, 431.
%\bibitem{ko} Kobayashi, Shoshichi, Differential geometry of complex vector bundles.
%Publications of the Mathematical Society of Japan, 15. Kano Memorial Lectures, 5. Princeton University Press, Princeton, NJ; Iwanami Shoten, Tokyo, 1987. xii+305 pp.
\bibitem{lo}
Lange, ~H., Ortega,  ~A., Prym varieties of cyclic coverings, 
 Geom. Dedicata 150  (2011), 391--403.

\bibitem{mn} Marcucci,  ~V, Naranjo,  ~J.C., Prym varieties of double coverings of elliptic curves,  Int. Math. Res. Notices 6 (2014), 1689-1698.


\bibitem{mp} Marcucci, ~V,  Pirola, ~G.P., Generic Torelli for Prym varieties of ramified coverings, Compos. Math. 148 (2012), 1147-1170.

\bibitem{nr} Nagaraj,~D.S.,  Ramanan, ~S., Polarisations of type $(1,2,\dots,2)$ on abelian varieties, Duke Math. J.  80 (1995), 157-194.   
\bibitem{no} Naranjo, ~J.C., Ortega, ~ A., Verra, ~A., Generic injectivity of the Prym map for double ramified coverings. Trans of AMS. DOI: 10.1090/tran/7459.  

\end{thebibliography}
\end{document}